\documentclass[12pt]{article}
\usepackage[cp866]{inputenc}
\usepackage{amssymb,amsmath,amsfonts}

\begin{document}

\pagestyle{plain}

\title{ On D. I. Moldavanskii's question about $p$-separable subgroups of a free group}

\author{Valerij. G. Bardakov}


\maketitle

\begin{abstract}

We prove that every nonabelian free group has a finitely generated isolated subgroup
which is not separable in
the class of nilpotent groups.

This enables us to give a negative  answer to the following
question by D.~I.~Moldavanskii in the ``Kourovka Notebook'': Is it
true that any finitely generated $p'$--iso\-la\-ted subgroup of a
free group is separable in the class of finite $p$--groups?

\end{abstract}


According to A. I. Mal'cev [1], a subgroup $H$
of a group $G$ is said to be {\it finitely separable from an element} $g \in G \setminus
H$ if there is a homomorphism
$\varphi$ of $G$ to some finite group such that
$\varphi(g) \not \in \varphi(H)$.

If a subgroup $H$ of a given group is finitely separable from every element
of this group which is not in $H$ we will say
that $H$ is  {\it finitely separable}. If we will consider
homomorphisms not on a finite group but on groups from some class
$\mathcal{K}$ we will define {\it separability in the class}
$\mathcal{K}$. In particular, if $\mathcal{K}$ is the class of
finite
$p$--groups then we will say that a subgroup $H$ is
$p$--{\it separable}. In this article $p$ is a fixed prime number.
The problem of finite separability is closely related to the generalized word problem (see [1]).
The {\it generalized word problem} for $H$ in $G$ asks for an algorithm
that decides whether or not the elements of $G$ lie in $H$.

M.~Hall proved in [2] that every finitely generated subgroup of a free  group
is finitely separable. One easily  sees that an analogue of Hall's theorem is not true
for the class of all finite $p$--groups. Indeed,
let $G = \langle a \rangle $ be an infinite cyclic group and $H = \langle a^q \rangle $ be a subgroup of
$G$, where $q \not = p$ is a prime number. Evidently, $a \not \in H$ but
for every homomorphism from $G$ onto a $p$--group the image of $a$ lies
in the image of $H$. Therefore, $H$ is not $p$--separable.

In the previous example $H$ is not $p'$--isolated in $G$.
Recall that a subgroup $H$ of $G$ is called $p'$--{\it
isolated}, if for every prime $q \not =p$ and for
every $g \in G$ an inclusion $g^q \in H$ implies that $g \in H$.
If in the previous example we will consider only
$p'$--isolated subgroups, then each of them will be
$p$--separable.

D.~I.~Moldavanskii suggested that the last fact is true
in every free group.

{\bf Problem} ([3, Problem 15.60]). {\it  Is it true that any finitely generated
$p'$--iso\-la\-ted
subgroup of a free group is separable in the class of finite $p$--groups?
It is easy to see that this is true for cyclic subgroups.}

E.~D.~Loginova [4, \S3] proves that in each finitely generated nilpotent
group every $p'$--isolated subgroup is $p$--separable. This
theorem says that the Moldavanskii's hypothesis may be true.

In the present article we prove the following

{\bf Theorem.} {\it In every free nonabelian group there is a
finitely generated isolated subgroup which is not separable in the
class of nilpotent groups.}

Recall that every finite $p$--group is nilpotent [5, p.
115] and every isolated subgroup is
$p'$--isolated for every prime  number $p$. Therefore, we obtain as consequence of the Theorem
that D.~I.~Moldavanskii's conjecture is not true.

{\bf Corollary.} {\it Let $p$ be a prime number. In every free
nonabelian  group there is a finitely generated $p'$--isolated
subgroup which is not separable in the class of finite
$p$--groups.}

In now we can start to prove the Theorem. Consider a free
nonabelian group
$$
F = \langle x, y, z_i (i \in I) \rangle,
$$
where $I$ is some index set (possibly empty). Let $H$ be a
subgroup of $G$:
$$
H = \langle x \, [y,\, x],~ y,~ z_j (j \in J) \rangle,
$$
where $[y, \, x] = y^{-1} \, x^{-1} \, y \, x$ and $J$ is a subset of $I$.
We see that $H$ is a proper subgroup of $F$ and, in particular, the element $x$ does not lie in $H$.
Indeed, let $a = x \, [y, \, x]$, $b = y$. Obviously,  every
nonempty word, reduced over the alphabet
$A = \{ a^{\pm 1},~ b^{\pm 1},~ z_j^{\pm 1} (j \in J) \}$, is reduced over the alphabet
$X = \{ x^{\pm 1},~ y^{\pm 1},~ z_i^{\pm 1} (i \in I) \}$ and
contains the letter $y$ or $y^{-1}$. Therefore, none of these words coincide with $x$.

Let us prove now that the subgroup $H$ is not separable in the class of nilpotent groups.
To do this we consider the quotient groups $F / \gamma_n F$, where
$$
\gamma_1 F = F, ~~~\gamma_{i+1} F = [\gamma_i F, F],~~~i = 1, 2,
\ldots ,
$$
are the terms of the lower central series of  $F$. We also consider
the natural homomorphisms
$$
\varphi_n : F \longrightarrow F/\gamma_n F,~~~ n = 1, 2, \ldots .
$$
Note that for each $\varphi_n$ the image $\varphi_n (H)$
is equal to the image of the subgroup
$\langle x,$ $y, z_j (j \in J) \rangle $.
The statement is obvious for $n = 1, 2$. For the case when $n > 2$ we apply the
following well-known fact [6, Theorem 31.2.5]:
If a nilpotent group $G$ is generated modulo
the commutator subgroup $G' = \gamma_2 G$
by elements $a_1, a_2, \ldots, a_r$,  then the elements $a_1, a_2, \ldots, a_r$
also generate $G$. Since all  $F
/ \gamma_n F$ are free nilpotent groups, it follows that $H$
is not separable in the class of all free nilpotent groups and hence $H$ is
not separable in the class of all nilpotent groups.

To complete the proof of the Theorem we must prove that $H$ is isolated in $F$.
For simplicity sake, we consider the case in which $F$ does not include the generators $z_i$.
The proof is similar in the general case.

{\bf Lemma.} {\it The subgroup $H = \langle x \, [y, \, x],~ y \rangle $ is isolated in
the free group}
$F = \langle x,~ y \rangle $.

{\bf Proof.} Let $a = x \, [y, \, x]$, $b = y$. We will consider
the elements of $F$ as words over the alphabet $X = \{
x^{\pm 1}, y^{\pm 1} \}$. The elements of $H$ will be considered both as
words over the alphabet $A = \{ a^{\pm 1}, b^{\pm 1} \}$ and the as words over the alphabet $X$
on substituting for $a$ and $b$ their expressions over the alphabet
$X$. By $``\equiv "$ will denote the graphical equality of words over the alphabet $X$, by $``="$
will denote equality in $F$.

Let for some reduced word $f$ from $F$ and for some integer number $m$ element $f^m$
lies in $H$. Obviously we may consider only the case $m > 0$. Then for
some integers $\alpha_1, \beta_1, \alpha_2, \beta_2, \ldots,
\alpha_k, \beta_k$, which are all nonzero except for, possibly,
$\alpha_1$ and $\beta_k$ we have that
$$
f^m = a^{\alpha_1} \, b^{\beta_1} \, a^{\alpha_2}  \, b^{\beta_2} \,
\ldots \,
a^{\alpha_k} \, b^{\beta_k}. \eqno{(1)}
$$
Note that if the word on the right--hand side of this equality is reduced over the alphabet $A$,
then it also reduced over the alphabet $X$. We will show that
 the reduced word $f^m$ might be considered as cyclically reduced.

Indeed, consider the last letter of
$f^m$. If it equals $y$, then $\beta_k > 0$. If the first letter of $f^m$
is not $y^{-1}$, then the word $f^m$ is cyclically reduced.
If the first letter of $f^m$ is $y^{-1}$ then
$\alpha_1 = 0$ and $\beta_1 < 0$. Therefore, the word on the right--hand side of (1)
is not cyclically reduced over $A$. Conjugate both sides of (1) by the element $y^{-1} =
b^{-1}$. We obtain an element in $H$ which is the $m$th power of an element
in
$F$. The case when the last letter of
$f^m$ is $y^{-1}$ is considered in a similar way.

Let $x$ be the last letter of
$f^m$.
Then $\beta_k = 0$, $\alpha_k > 0$. If the first letter of
$f^m$ is not $x^{-1}$, then $f^m$ is cyclically reduced. If the first letter of $f^m$
is $x^{-1}$, then $\alpha_1
< 0$. Conjugating both sides of (1) by $a^{-1}$, we again find an
element in $H$ which is the $m$th power of some element in $F$.
In the case in which the last letter of $f^m$ is $x^{-1}$ the argument is similar.
Repeating this process, we will arrive to an equality of the same form as (1)
where the word on the left--hand side is cyclically reduced, as desired.

So, we are assuming from now on that the word $f^m$ on the left--hand side of (1)
is cyclically reduced.
In this case $f^m \equiv f_1 \, f_2 \, \ldots \, f_m$,
where $f_i \equiv f$, i.e., there are no reductions in the products $f_i \, f_{i+1}$,
$i = 1, 2, \ldots, m-1$. Consider the equality
$$
f_1 \, f_2 \, \ldots  \, f_m = a^{\alpha_1} \, b^{\beta_1} \,
a^{\alpha_2} \,
b^{\beta_2} \, \ldots \, a^{\alpha_k} \, b^{\beta_k}. \eqno{(2)}
$$
If $f_1$ is equal to an initial subword of the word from the right--hand side of (2)
then $f_1$ is equal to  some word over the alphabet
 $A$, i.e., $f = f_1 \in H$ and the statement is true.

We will prove that other possibilities are impossible.

Suppose that $f_1$ is a product of an initial subword on the right--hand side of (2)
and some subword of $a$ or $a^{-1}$, i.e., for some exponents $\alpha_l$, $1 \leq l
\leq k,$ we have $a^{\alpha_l} \equiv a^{\gamma} \, a_0 \, a_1 \,
a^{\delta}$, $|\alpha_l| = |\gamma | + |\delta | + 1,$ and
$$
f_1 \equiv a^{\alpha_1} \, b^{\beta_1} \, a^{\alpha_2} \,
b^{\beta_2} \,
\ldots \, b^{\beta_{l-1}} \, a^{\gamma } \, a_0,
$$
where $a \equiv a_0 \, a_1$ if $\alpha_l > 0$ and $a^{-1} \equiv
a_0 \,
a_1$ if $\alpha_l < 0$. The first letter of
$f_1$ may be one of the set $\{ x^{\pm 1},~~ y^{\pm 1} \}$.
Consider all these possibilities.

Let the first letter of $f_1$ be $x$. Since $f_1 \equiv
f_2$ it follows that $f_2$ must begin with $x$ as well and, moreover,
either
$$
f_1 \equiv a^{\alpha_1} \, b^{\beta_1} \, a^{\alpha_2} \, b^{\beta_2}
\, \ldots \, b^{\beta_{l-1}} \, a^{\gamma } \, x \, y^{-1} \, x^{-1} \, y,~~~\gamma
\geq 0,~~~a_0 \equiv x \, y^{-1} \, x^{-1} \, y,
$$
or
$$
f_1 \equiv a^{\alpha_1} \, b^{\beta_1} \, a^{\alpha_2} \, b^{\beta_2}
\, \ldots \, b^{\beta_{l-1}} \, a^{\gamma } \, x^{-1} \, y^{-1},~~~\gamma \leq
0,~~~a_0 \equiv x^{-1} \, y^{-1}.
$$
In the first case we have
$$
f_2 \equiv x \, a^{\delta} \, b^{\beta_l} \, a^{\alpha_{l+1}} \, w,
$$
where $w$ is some reduced word. Since the first letter of $f_1$ is $x$, we have $\alpha_1 > 0$
and
$$
f_1 \equiv x \, y^{-1} \, x^{-1} \, y \, x \, a^{\alpha_1 - 1} \, b^{\beta_1} \, \ldots
\, a^{\gamma} \, x \, y^{-1} \, x^{-1} \, y.
$$
In order to the first three letters of $f_2$ are equal to the first three letters of $f_1$
it is necessary that $\delta = 0$, $\beta_l
= -1$, $\alpha_{l+1} < 0$. But in this case
$$
f_2 \equiv x \, y^{-1} \, x^{-1} \, y^{-1} \, x \, y \, x^{-1} \, a^{\alpha_{l+1} + 1}
\, w.
$$
By comparing the fourth letters of $f_1$ and $f_2$, we see that $f_1 \not
\equiv f_2$. Therefore this case is impossible.

The case in which the first letter of $f_1$ is $x^{-1}$ considered similarly.

Let the first letter of  $f_1$ be $y$. Then $\alpha_1 = 0$, $\beta_1 > 0$.
Since the first letter of $f_2$ is $y,$ too; therefore, either
$$
f_1 \equiv b^{\beta_1} \, a^{\alpha_2} \, b^{\beta_2} \, \ldots \, a^{\gamma}
\, x \, y^{-1} \, x^{-1},~~~\gamma \geq 0,
$$
or
$$
f_1 \equiv b^{\beta_1} \, a^{\alpha_2} \, b^{\beta_2} \, \ldots \, a^{\gamma}
\, x^{-1} \, y^{-1} \, x,~~~\gamma \leq 0.
$$
Clearly, in both cases $\beta_1 = 1$. Hence, in the
first case
$$
f_2 \equiv y \, x \, a^{\delta} \, b^{\beta_l} \, a^{\alpha_{l+1}} \, w,~~~\delta
\geq 0,
$$
where $w$ is some reduced word. In order that the second letter of $f_1$
is equal to the second letter of $f_2$ it is necessary that
$\alpha_2 > 0$, i.e.,
$$
f_1 \equiv y \, x \, y^{-1} \, x^{-1} \, y \, x \, a^{\alpha_2-1} \, b^{\beta_2} \, \ldots
\, a^{\gamma} \, x \, y^{-1} \, x^{-1}.
$$
In order that the third and fourth letters in the word $f_2$ are equal to the third and fourth
letters of $f_1$ accordingly, it is necessary that $\delta =
0$, $\beta_l = -1$, $\alpha_{l+1} < 0$, i.e.,
$$
f_2 \equiv y \, x \, y^{-1} \, x^{-1} \, y^{-1} \, x \, y \, x \, a^{\alpha_{l+1}-1} \, w,
$$
but if we compare the fifth letters  $f_1$ and $f_2$ then we will see that $f_1 \not
\equiv f_2$. In the second case
$$
f_2 \equiv y \, x^{-1} \, a^{\delta} \, b^{\beta_l} \, a^{\alpha_{l+1}}
\, w,~~~\delta \leq 0.
$$
In order that the second letter of $f_1$ is equal to the second letter of $f_2$
it is necessary that $\alpha_2 < 0$, i.e.,
$$
f_1 \equiv y \, x^{-1} \, y^{-1} \, x \, y \, x^{-1} \, a^{\alpha_2+1} \, b^{\beta_2}
\, \ldots \, a^{\gamma} \, x^{-1} \, y^{-1} \, x.
$$
In order that the third and fourth letters of $f_2$
are equal to the third and fourth letters of $f_1$
accordingly; it is necessary that
 $\delta =
0$, $\beta_l = -1$, $\alpha_{l+1} > 0$, i.e.,
$$
f_2 \equiv y \, x^{-1} \, y^{-1} \, x \, y^{-1} \, x^{-1} \, y \, x \, a^{\alpha_{l+1}-1}
\, w,
$$
but if we compare the fifth letters of $f_1$ and  $f_2$ then we will see that $f_1 \not
\equiv f_2$. Therefore $f_1$ cannot begin with $y$.

Arguing similarly, we see that it is not possible for $f_1$ to begin with $y^{-1}$.
This completes the proof of the lemma and  theorem.

{\bf {\it Acknowledgements.}}
I am very grateful to M.~V.~Neshchadim for very stimulating
discussions. Special thanks go to the participants of the
seminar ``Evariste Galois'' at Novosibirsk State
University for their kind attention to my work.

\vskip 24pt

 \vskip 20pt

\centerline{\bf REFERENCES} \vskip 12pt
\begin{enumerate}
\item
    A. I. Mal'cev, On homomorphisms onto finite groups, Uchen. Zapiski Ivanovsk. Ped. instituta,
18, N~5 (1958), 49--60 (also in ``Selected Papers'', Vol. 1, Algebra, 1976, 450--462) (Russian).
\item
    M.~Hall, Jr., Coset representations in free groups, Trans. Amer. Math. Soc.,
67, N~2 (1949), 421--432.
\item
   The Kourovka Notebook (Unsolved Problems
in Group Theory), 15th ed., Institute of Mathematics SO RAN, Novosibirsk, 2002.
\item
E.~D.~Loginova, Residual finiteness of the free product of two groups with commuting subgroups.
 Sibirsk. Mat. Zh., 40, N~2, (1999), 395-407 (Russian).
\item
        M.~I.~Kargapolov, Yu.~I.~Merzljakov, Fundamentals of the Theory of Groups, New York:
Springer, 1979.
\item
    H. Neumann, Varieties of Groups, Ergebnisse der Mathematik und ihrer Grenzgebiete, Band 37,
Springer--Verlag, Berlin--Heidelberg--New York, 1967.

\end{enumerate}

\bigskip
\bigskip
\noindent
Author address:

\bigskip

\noindent
Valerij G. Bardakov\\ Sobolev Institute of Mathematics,\\
Novosibirsk, 630090, Russia\\ {\tt bardakov@math.nsc.ru}

\end{document}